\documentclass[11pt,reqno,tbtags]{amsart}
\usepackage{amssymb}
\usepackage{natbib}
\bibpunct[, ]{[}{]}{;}{n}{,}{,}

\allowdisplaybreaks

\newtheorem{theorem}{Theorem}

\newtheorem{lemma}[theorem]{Lemma}

\theoremstyle{definition}

\theoremstyle{remark}

\newcounter{thmenumerate}

\newcounter{xenumerate}   

\newcommand\marginal[1]{\marginpar{\raggedright\parindent=0pt\tiny #1}}

\begingroup
  \divide\count255 by 60
  \multiply\count255 by -60
  \advance\count255 by \time
  \ifnum \count255 < 10 \xdef\klockan{\the\count1.0\the\count255}
  \else\xdef\klockan{\the\count1.\the\count255}\fi
\endgroup

\def\rompar(#1){\textup(#1\textup)}    

\def\xexp(#1){e^{#1}}

\newcommand\ntoo{\ensuremath{{n\to\infty}}}

\newcommand\eg{e.g.\spacefactor=1000}

\newcommand\whp{{whp}}

\newcommand\gd{\delta}

\newcommand\gG{\Gamma}

\newcommand\eps{\varepsilon}

\newcommand\cC{\mathcal C}

\newcommand\cF{\mathcal F}

\def\[#1]{[\![#1]\!]}

\renewcommand{\=}{:=}

\newcommand{\gnc}{{\mathcal G}(n, c/n)}
\newcommand{\gnd}{{\mathcal G}_d(n)}

\newcommand\REM[1]{\texttt{[#1]}\marginal{XXX}}

\newcommand\urladdrx[1]{{\urladdr{\def~{{\tiny$\sim$}}#1}}}

\begin{document}
\title
{Dismantling sparse random graphs}

\date{September 12, 2007} 

\author{Svante Janson}
\address{Department of Mathematics, Uppsala University, PO Box 480,
SE-751~06 Uppsala, Sweden}
\email{svante.janson@math.uu.se}
\urladdrx{http://www.math.uu.se/~svante/}
\author{Andrew Thomason}
\address{DPMMS, Centre for Mathematical Sciences, Wilberforce Road,
Cambridge CB3 0WB, United Kingdom}
\email{a.g.thomason@dpmms.cam.ac.uk}


\begin{abstract} 
We consider the number of vertices that must be removed from a graph $G$ in
order that the remaining subgraph has no component with more than $k$ vertices.
Our principal observation is that, if $G$ is a sparse random graph or a random
regular graph on $n$ vertices with $n\to\infty$, then the number in question is
essentially the same for all values of $k$ that satisfy both $k\to\infty$ and
$k=o(n)$.
\end{abstract}

\maketitle

The process of removing vertices from a graph $G$ so that the remaining
subgraph has only small components is known as \emph{fragmentation}. Typically,
the aim is to remove the least possible number of vertices to achieve a given
component size; this is equivalent to determining the largest induced subgraph
whose components are at most that size. This process has been studied in (at
least) two different lines of research, from different perspectives 
and with quite different component sizes. In this
note we point out that, as far as sparse random graphs are concerned, these
two perspectives actually arrive at the same answer.

Let $\Gamma$ be a class of graphs. The classes we shall mostly be interested in
are the classes ${\mathcal C}_k$, the class of graphs whose components have at
most $k$~vertices, and ${\mathcal F}$, the class of forests.  Given such a
class $\Gamma$, we define
\begin{equation*}
N(G,\Gamma)\=\max\{|S|\,:\,G[S]\in\Gamma\},
\end{equation*}
where $S$ is a subset of the vertices of $G$ and $G[S]$ denotes the subgraph of
$G$ induced by $S$. We also define
\begin{equation*}
\nu(G,\Gamma)\=N(G,\Gamma)/|G|,
\end{equation*}
so that $0\le \nu(G,\Gamma)\le 1$. 
(To make this always defined, we set $N(G,\gG)=0$ if no induced
subgraph of $G$ belongs to $\gG$; equivalently, we may regard the
empty graph with no vertices as an element of $\gG$.) 
Thus, for example, the size
of a largest independent set in $G$ is 
$N(G, {\mathcal C}_1)=\nu(G, {\mathcal C}_1)|G|$.
(This is known as the \emph{independence number}.)
Similarly, 
$n-N(G, \cF)$ is the \emph{decycling number}, see \eg{} \citet{Karp}.

In this notation, the study of fragmentation is the study of the parameter
$\nu(G, {\mathcal C}_k)$ for various values of~$k$. From the point of view of
graph theory, it is natural to consider $\nu(G, {\mathcal C}_k)$ for some large
but finite value of~$k$, for graphs $G$ in which the number of vertices $n=|G|$
grows large. This study was initiated by Edwards and Farr~\cite{ef01,ef04}. On
the other hand, in the study of vaccination, see for example
\citet{SJ199} and the references therein,
the vertices of the graph are individuals in some
population, with edges representing the opportunity of passing on a
disease. If a vertex is vaccinated it becomes unable to spread the disease; a
vaccination strategy is a way to ensure that the subgraph induced by the
unvaccinated vertices has only small components (relative to the total
population). The vaccinator is thus interested in $\nu(G, {\mathcal C}_{\delta
n})$ for small values of~$\delta$. 
(For further details and for variations on this theme, see
\citep{SJ199}.)

In both studies it is natural to consider the behaviour of these parameters on
two standard models of sparse random graphs. Let ${\mathcal G}(n, c/n)$ denote
the probability space of graphs with vertex set $\{1,\ldots,n\}$ with edges
chosen independently with probability $c/n$, and let ${\mathcal G}_d(n)$ denote
the space of $d$-regular graphs on the same vertex set. We shall assume that
$c>1$ and $d\ge3$ are fixed as $n\to\infty$.  
(For odd $d$, we of course have to assume that $n$ is even.)
The main observation of this
paper is that, for graphs in these spaces, the graph theoretic approach and the
vaccination approach arrive at the same answer: that is, perhaps surprisingly,
fragmenting into large but finite components \whp{} costs no more than
just fragmenting into components of size~$o(|G|)$. 
(A sequence of events $(A_n)$ is said to hold \whp{} if $\Pr (A_n)\to1$ as 
$n\to\infty$.)

\begin{theorem}\label{thetheorem}
Let $c>1$, $d\ge3$ and $\epsilon>0$ be given. Then there exists $\delta>0$ such
that, if $G\in {\mathcal G}(n, c/n)$ or $G\in{\mathcal G}_d(n)$ then
\begin{equation*}
\nu(G, {\mathcal C}_{\delta n})\,\le\,\nu(G, {\mathcal C}_{1/\delta}) +
\epsilon
\end{equation*}
holds \whp{} as $n\to\infty$.
\end{theorem}

Before giving the proof of this theorem, we make a few more remarks.
Edwards and Farr~\cite{ef01,ef04} considered general graphs of bounded maximum
degree; in particular they studied the parameter
\begin{equation*}
\beta_d\,\=\,\sup_k \min\{\nu(G, {\mathcal C}_k)\,:\, G\mbox{ has maximum
  degree } d\,\}
\end{equation*}
(note that the $\alpha_d$ of Edwards and Farr equals $1-\beta_d$). One way to
think of this parameter is that, if $\beta<\beta_d$, then there is some finite
$k$ for which every graph $G$ of maximum degree~$d$ has an induced subgraph
$G[S]$ with at least $\beta |G|$ vertices but with no component larger
than~$k$. Trivially $\beta_1=\beta_2=1$. and it is shown in~\cite{ef01} that
$\beta_3=\frac{3}{4}$. In general they showed that $\beta_d\ge\frac{3}{d+1}$;
a complementary inequality  $\beta_{2d}\le\frac{2}{d+1}$ was
proved by Haxell, Pikhurko and Thomason~\cite{hpt} (so answering
affirmatively the question posed by Edwards and Farr~\cite{ef03}
as to whether $\beta_d\to0$ as $d\to\infty$).

The parameter $\nu(G,{\mathcal F})$, describing the largest induced forest in
$G$, is significant in the study of fragmentation, because any forest $F$ is
easily fragmented by removing a few vertices. The following simple lemma is
given only because it is best possible, as exemplified by a path.
\begin{lemma}\label{thelemma}
If $F$ is a forest then $\nu(F,{\mathcal C}_k)\ge 1- (k+1)^{-1}$.
\end{lemma}
\begin{proof}
We may assume that $F$ is a tree, and proceed by induction on $n=|F|$, the case
$n\le k+1$ being trivial. For larger~$n$, note that the removal of any edge
leaves two components. Orient the edge towards the larger of these (break ties
arbitrarily) and colour the edge red if both components have more than~$k$
vertices. If there are no red edges, remove a sink vertex (there must be one
since $F$ is acyclic) and observe that 
this leaves only components with at most $k$ vertices, and thus
$\nu(F,{\mathcal C}_k)=1-1/n \ge 1-1/(k+1)$. 
If there are
red edges, it is easy to see that they form a connected subgraph and so a tree;
the removal of a leaf vertex of the red tree breaks $F$ into a tree with at
most $n-k-1$ vertices plus some components of size at most~$k$, and the proof
follows by induction.
\end{proof}
The parameter $\nu(G,{\mathcal P})$, where $\mathcal P$ is the class of planar
graphs, is likewise significant for fragmentation, because a planar graph $P$
can be fragmented quite efficiently by means of the separator theorem of Lipton
and Tarjan~\cite{liptonT}, and in fact $\nu(P,{\mathcal C}_k)\ge 1-24k^{-1/2}$
holds (see~\cite{ef01} or~\cite{edwards+mcdiarmid:94}). From this and from
Lemma~\ref{thelemma} it follows that, for any graph~$G$, both $\nu(G,{\mathcal
C}_k)\ge\nu(G,{\mathcal F}) - (k+1)^{-1}$ and $\nu(G,{\mathcal
C}_k)\ge\nu(G,{\mathcal P}) - 24 k^{-1/2}$ hold. On the other hand, if $n_k(G)$
denotes the number of cyles in $G$ of length at most~$k$, then (by removing a
vertex from each cycle) we have $\nu(G,{\mathcal P})\ge\nu(G,{\mathcal F})\ge
\nu(G,{\mathcal C}_k)-n_k(G)/|G|$. Thus if we restrict our attention to large
graphs $G$ for which $n_k(G)=o(|G|)$ for each fixed $k$,
then the three parameters $\nu(G,{\mathcal
F})$, $\nu(G,{\mathcal P})$ and $\nu(G,{\mathcal C}_k)$ are asymptotically the
same for large~$k$. Graphs in ${\mathcal G}(n, c/n)$ or ${\mathcal G}_d(n)$
enjoy this property \whp{}
(see \cite{bollobas} or \cite{JLR}).

The parameters $\nu(G,{\mathcal C}_1)$ and $\nu(G,{\mathcal F})$ for  $G\in
{\mathcal G}(n, c/n)$ and $G\in{\mathcal G}_d(n)$ have already
received considerable attention. The first of these (the 
independence number) was studied by
\citet{f90} for ${\mathcal G}(n, c/n)$, 
see also \cite[Section 7.1]{JLR},
and 
Frieze and \L{}uczak~\cite{fl} for
$G\in{\mathcal G}_d(n)$, 
and
information on the second is given by Bau, Wormald and Zhou~\cite{bwz}.  In
fact it is shown in~\cite{fl} that, for $G\in{\mathcal G}_d(n)$,
$\nu(G,{\mathcal C}_1)\sim 2\log d / d$ holds \whp{} for large~$d$,
which is already enough to answer the above-mentioned question about the limit
of $\beta_d$, and in~\cite{hpt} it is verified that $\nu(G,{\mathcal F})\sim
2\log d/d$ \whp.
(These statements involve double limits, as first $n\to\infty$ and
then $d\to\infty$. More precisely, by 
\cite{fl,hpt},
for every $\eps>0$, there exists
$d_\eps$ such that for $\gnd$ with any fixed $d\ge d_\eps$, \whp{} holds
$(2-\eps)\log d / d
\le
\nu(G,{\mathcal C}_1)
\le
\nu(G,\cF)
\le (2+\eps)\log d / d$ 
as \ntoo. A similar result holds for $\gnc$, by
\cite{f90}
and a first moment argument as in \cite{hpt}.)

We are now ready for the proof of Theorem~\ref{thetheorem}.

\begin{proof}[Proof of Theorem~\ref{thetheorem}]
We claim that the following holds \whp, if $\delta$ is small
enough: 
\emph{Each set $T$ of at most $\delta n$ vertices spans at most
$(1+\epsilon/3)|T|$ edges}.

The theorem follows from this claim; for let $S$ be a set such that
$G[S]\in{\mathcal C}_{\delta n}$ and $|S|=\nu(G, {\mathcal C}_{\delta n})$. By
the claim, from each component $G[T]$ of $G[S]$ we may remove at most
$\epsilon|T|/3$ edges so that it becomes acyclic or unicyclic; thus, by
removing $\epsilon|S|/3$ edges we can make all components of $G[S]$ acyclic or
unicylic. There are at most $\epsilon|S|/3$ components of size larger than
$3/\epsilon$, and so by removing a further $\epsilon|S|/3$ edges we can make
all these large components acyclic. 
Hence, removing vertices instead of edges,
there exists $S^\prime\subset S$,
$|S^\prime|\ge|S|-2\epsilon n/3$, such that $G[S^\prime]$ consists of a forest
$F$ plus components of size at most~$3/\epsilon$. By Lemma~\ref{thelemma},
the removal of a further $\epsilon|S^\prime|/3$ vertices from $G[S^\prime]$
leaves only components of size at most~$3/\epsilon$. Therefore $\nu(G,
{\mathcal C}_{3/\epsilon}) \ge \nu(G, {\mathcal C}_{\delta n})-\epsilon$. We
can of course assume that $\delta<\epsilon/3$, and this proves the theorem.

To prove the claim, consider first the case $G\in {\mathcal G}(n, c/n)$. Let
$T$ be a set of
$t\le n\gd$ vertices, and let $\tau=n/t\ge 1/\delta$; we can make
$\tau$ large by making $\delta$ small. Let $X$ be the random variable counting
the number of edges in $G[T]$. Then $X$ is binomially distributed with mean
$\lambda=\frac{c}{n}\binom{t}{2}\le
\frac{ct}{2\tau}$. By the version of the Chernoff bound in
\cite[Corollary~2.4]{JLR}, if $x=m\lambda>\lambda$, then
${\mathbb P}(X\ge x)\le\exp(-lx)$, where $l=\log m-1+1/m$. Taking
$x=(1+\epsilon/4)t$ we have $m=x/\lambda\ge(1+\epsilon/4)2\tau/c>1$ if
$\delta<2/c$, and 
$l>\log m-1\ge\log \tau-1-\log c$. The number of sets of size $t$ is
$\binom{n}{t}\le(e\tau)^t$. Therefore the probability $P_t$ that the claim
fails for some set $T$ of size~$t$ satisfies
\begin{equation*}
P_t\,\le\,\exp\{t(1+\log\tau)-t(1+\epsilon/4)(\log \tau - 1 - \log c)\}
\,\le\,\exp\{-\epsilon t(\log \tau)/8\}
\end{equation*}
if $\delta$ is small enough. For $t\ge\log n$ we have $P_t\le n^{-2}$ if
$\delta$ is small, and for $t\le \log n$ we have
$\log\tau\ge(\log n)/2$ so 
$P_t\le n^{-\epsilon/16}$. Thus $\sum_{1\le t\le \delta n}P_t=o(1)$, which
proves the claim.

The case $G\in{\mathcal G}_d(n)$ is similar. The calculation is messier but
fortunately we need not give it here, because it is essentially that of the
proof of Lemma~5.1 of Janson and Luczak~\cite{jl}. They prove (see Remark~5.2)
that each set $T$ with $|T|\le\delta n$ has average degree less than~$k$, where
$k\ge3$. However their interest was in integer values of~$k$, and the proof,
in which $k$ appears everywhere as a variable, works perfectly well for any
fixed $k>2$, which is exactly our claim.
\end{proof}

In the light of Theorem~\ref{thetheorem} it is interesting to consider the
fragmentation of random graphs $G$ into components of size proportional to
$|G|$. Now, given $c>0$ and $0<x\le 1$, the inequalities of
Azuma--Hoeffding~\cite[Corollary 2.27]{JLR} or of Talagrand~\cite[Theorem 2.29
and Remark 2.36]{JLR} show that the value of $\nu(G, {\mathcal C}_{x n})$ for
$G\in {\mathcal G}(n, c/n)$ is highly concentrated. It seems very likely that
there is a real number $f_c(x)$ such that $\nu(G, {\mathcal C}_{x n})$ tends to
$f_c(x)$ in probability; that is, for all $\epsilon>0$, $\Pr\{|\nu(G, {\mathcal
C}_{x n}) - f_c(x)| > \epsilon\}\to0$ as $n\to\infty$.  But it is not actually
known whether $f_c(x)$ exists; this is an unfortunate state of affairs shared
with many standard parameters such as $\nu(G,{\mathcal C}_1)$ and
$\nu(G,{\mathcal F})$. 
(An exception here is $\nu(G,{\mathcal C}_1)$ when $c\le e$; 
see~\cite[Corollary~1]{GNS}.)
Nevertheless, 
our final comments can be stated more cleanly
by assuming both that $f_c(x)$ exists, and also that, for each fixed~$k$,
$\nu(G, {\mathcal C}_k)$ tends to a limit $\phi_c(k)$ in probability. 
Corresponding limits will be assumed for $G\in{\mathcal G}_d(n)$ too,
and we denote these by $g_d(x)$ and $\gamma_d(k)$.
(Note that, at least, one can show that
$\nu(G, {\mathcal C}_{x n})$ and
$\nu(G, {\mathcal C}_{k})$ are highly concentrated
for $G\in \gnd$ too, using
a version for
random permutations of the inequality by Talagrand
\cite[Theorem 5.1]{Talagrand95}, see also \citet[Theorem 1.1]{McD02}, and
arguing as in \citet{hpt}.)

We interpret $f_c(x)$ as meaning that the largest induced subgraph of a random
graph $G\in {\mathcal G}(n, c/n)$ having no component larger than $xn$
has
(about)
$f_c(x)n$ vertices, and the largest induced subgraph having no component larger
than~$k$ has $\phi_c(k)n$ vertices. The numbers $g_d(x)$ and $\gamma_d(k)$ have
similar interpretations. All these functions are increasing in their
arguments. So the limits $\lim_{x\to0}f_c(x)$ and $\lim_{x\to0}g_d(x)$ exist;
we define $f_c(0)$ and $g_d(0)$ to be the values of these limits, thus making
$f_c$ and $g_d$ continuous at zero. The content of Theorem~\ref{thetheorem}
is then that
\begin{equation*}
\lim_{k\to\infty}\phi_c(k)=f_c(0)\=\lim_{x\to0}f_c(x)
\quad\mbox{and}\quad
\lim_{k\to\infty}\gamma_d(k)=g_d(0)\=\lim_{x\to0}g_d(x)\,.
\end{equation*}
[Note that $\lim\gamma_d(k)$ corresponds very closely to the parameter
$\beta_d$ defined earlier, the only difference being that $\beta_d$ takes
account of all $d$-regular graphs, whereas here we consider only almost all.]

The function $f_c$ satisfies the Lipschitz condition $f_c(x) - f_c(y)
<\frac{1}{y}(x-y)$ for $y<x$, since, if $G[S]$ has component sizes at most
$xn$, then at most $1/y$ of these are larger than $yn$, and each of these can
be reduced to size $yn$ by removing at most $(x-y)n$ vertices. This, together
with the fact that $f_c$ is increasing, shows that $f_c$ is a continuous
function on the unit interval, and similarly so is~$g_d$.

The famous theorem of Erd\H{o}s and R\'enyi~\cite{er}, that $G\in {\mathcal
 G}(n, c/n)$ almost certainly has a unique giant component of size 
$(1+o(1))\rho(c)n$, where $\rho(c)=1-e^{-c\rho(c)}$, means that $f_c(x)=1$ for
$\rho(c)\le x\le1$. There is no corresponding fact for random regular graphs,
of course; we just have $g_d(1)=1$.

It seems likely that the function $f_c$ is \emph{strictly} increasing on the
interval $[0,\rho(c)]$ and that $g_d$ is strictly increasing on $[0,1]$. This
would mean that continuous inverse functions $f^{-1}_c$ and $g_d^{-1}$
exist. For the vaccinator, the function $f^{-1}_c$ would be more natural than
$f_c$ itself; $f^{-1}_c(z)n$ is the smallest component size achievable by
vaccinating $(1-z)n$ people. But we cannot show strict monotonicity except at
the right-hand end of the range. 
It was proved by \citet[Theorem 3.9]{SJ178} that
for every $\epsilon>0$ there exists $\delta>0$ such that \whp{}
after removal of any $\gd n$ vertices from $\gnc$, there is still a
giant component of order at least $(\rho(c)-\eps/2)n$, and thus
$f_c(\rho(c)-\epsilon) < 1-\gd=f(\rho(c))-\delta$;
this is an easy consequence of the corresponding result for edge deletions
by
Luczak and McDiarmid~\cite[Lemma~2]{lm}, who
proved that for every $\epsilon>0$ there exists $\delta>0$ such that the giant
component \whp{} has no two sets, each of size at least $\epsilon n$,
having at most $\delta n$ edges between the sets. 
A similar argument can be given to
show that $g_d(1-\epsilon)< g_d(1)-\delta$.

In conclusion, our main open question is whether $f_c$ and $g_d$ are strictly
increasing. We would also like to know more about the subgraphs $G[S]$ of order
$f_c(x)n$ or $g_d(x)n$ that have no component larger than~$xn$: how many
components do they have? The reader who is interested in these questions can
readily formulate them in a way that does not involve the uncertain existence
of $f_c$ and $g_d$.

We finally remark that corresponding questions can be formulated for
removal of edges instead of vertices. In that case, the central
parameter is the largest number of edges in a (not necessarily
induced) subgraph of $G$ that belongs to $\gG$. For $\cF$, the class
of forests, this is easy (unlike the vertex case treated above, see
\citet{Karp}), but we see no easy answers for \eg{} $\cC_{\gd n}$ and
leave these versions as problems for the interested reader.

\section*{Acknowledgement}
This research was mainly done during a visit by SJ to the University
of Cambridge, partly funded by Trinity College.

\newcommand\AAP{\emph{Adv. Appl. Probab.} }
\newcommand\JAP{\emph{J. Appl. Probab.} }
\newcommand\JAMS{\emph{J. \AMS} }
\newcommand\MAMS{\emph{Memoirs \AMS} }
\newcommand\PAMS{\emph{Proc. \AMS} }
\newcommand\TAMS{\emph{Trans. \AMS} }
\newcommand\AnnMS{\emph{Ann. Math. Statist.} }
\newcommand\AnnPr{\emph{Ann. Probab.} }
\newcommand\CPC{\emph{Combin. Probab. Comput.} }
\newcommand\JMAA{\emph{J. Math. Anal. Appl.} }
\newcommand\RSA{\emph{Random Struct. Alg.} }
\newcommand\ZW{\emph{Z. Wahrsch. Verw. Gebiete} }
\newcommand\DMTCS{\jour{Discr. Math. Theor. Comput. Sci.} }

\newcommand\AMS{Amer. Math. Soc.}
\newcommand\Springer{Springer-Verlag}
\newcommand\Wiley{Wiley}

\newcommand\vol{\textbf}
\newcommand\jour{\emph}
\newcommand\book{\emph}
\newcommand\inbook{\emph}
\def\no#1#2,{\unskip#2, no. #1,} 

\newcommand\webcite[1]{\hfil\penalty0\texttt{\def~{\~{}}#1}\hfill\hfill}
\newcommand\webcitesvante{\webcite{http://www.math.uu.se/\~{}svante/papers/}}
\newcommand\arxiv[1]{\texttt{arXiv:#1}}

\def\nobibitem#1\par{}

\end{document}